\documentclass[11pt]{amsart}

\usepackage{amssymb, mathrsfs, mathabx, hyperref, geometry}
\usepackage[english,francais]{babel}

\newtheorem{theorem}{Theorem}[section]
\newtheorem{lemma}[theorem]{Lemma}
\newtheorem{e-proposition}[theorem]{Proposition}

\newtheorem{e-definition}[theorem]{Definition\rm}

\newcommand{\R}{\mathbb{R}}

\DeclareMathOperator{\Tr}{Tr}

\DeclareMathOperator{\vol}{vol}

\DeclareMathOperator{\Op}{Op}

\title{Quantum Ergodicity for Eisenstein functions}
\author{Yannick Bonthonneau}
\email{yannick.bonthonneau@gmail.com}
\address{CIRGET, UQ\`AM, 201 av. Pr\'esident Kennedy, Montr\'eal, Qu\'ebec, Canada, H2X 3Y7}
\author{Steve Zelditch}
\email{zelditch@math.northwestern.edu}
\address{Department of Mathematics, Northwestern  University, Evanston, IL 60208, USA}
\thanks{Research partially supported by NSF grant DMS-1541126}
\begin{document}

\begin{abstract}
\selectlanguage{english}
A new proof is given of Quantum Ergodicity for Eisenstein Series for cusped hyperbolic surfaces. This result is also extended to higher dimensional examples, with variable curvature.

\vskip 0.5\baselineskip

\selectlanguage{francais}
\noindent{\bf R\'esum\'e} \vskip 0.5\baselineskip \noindent
{\bf Ergodicit\'e Quantique des Fonctions d'Eisenstein. }
On donne une nouvelle preuve de l'Ergodicit\'e Quantique des s\'eries d'Eisenstein pour les surfaces de Riemann \`a pointes. On \'etend aussi ce r\'esultat en plus grande dimension, en autorisant la courbure variable.

\end{abstract}
\maketitle

\selectlanguage{english}

\section{Introduction}


The purpose of this note is to prove a quantum ergodicity theorem for Eisenstein
series on a rather general cusp manifold with ergodic geodesic flow. In the special
case of finite area hyperbolic surfaces with cusps, a similar quantum ergodicity theorem
was proved in  \cite{Zelditch-91}. In this note, we give a new, simpler  and more general proof of that result which generalizes to any cusped manifolds with ergodic geodesic
flow. In particular, the manifolds may  have variable negative curvature. 

To state the result, we introduce some notation and terminology.
A  $(d+1)$-dimensional \emph{cusp manifold} is a Riemannian manifold $M$ which decomposes as a compact manifold $M_0$ whose boundaries are torii, and a finite number $\kappa$ of topological half-cylinders $Z_1$, $\dots$, $Z_\kappa$, such that the metric on these $Z_i$'s is
\begin{equation}
ds^2 = \frac{ dy^2 + d\theta^2}{y^2}
\end{equation}
where $y$ is the vertical coordinate, starting at $a_i >0$, and $\theta$ is the horizontal coordinate, living in $\R^d/\Lambda_i$, where $\Lambda_i$ is a unimodular lattice. On such a manifold, we can pick a global coordinate $y_M$, that coincides with the vertical coordinate sufficiently high in the cusps, and is some constant in $M_0$.

For such a manifold, the spectrum of the Laplace-Beltrami operator $-\Delta$ acting on $L^2(M,g)$ decomposes into the pure spectrum part ($\lambda_0 = 0 < \lambda_1 \leq \dots$) possibly finite, possibly infinite, and the absolutely continuous spectrum $[1/4, + \infty )$.

To each eigenvalue $\lambda$, we can associate an eigenfunction $u_\lambda$ (repeating eventual multiple eigenvalues). For the continuous spectrum, we recall the existence of meromorphic families \cite{CdV-81,Muller-86} of smooth, non-$L^2$ eigenfunctions $E_i(s)$ that satisfy the following. For each $i$, $E_i$ decomposes as
\begin{equation}
E_i(s) = I_i y^s  + G_i(s)
\end{equation}
where $I_i$ is the characteristic function of $Z_i$, and $G_i(s)$ is a function in $L^2(M)$ whenever $\Re s > d/2$. We also require that for all $s$,
\begin{equation}
(-\Delta - s(d-s))E_i(s) = 0
\end{equation}
Such a collection is unique, and they are called the Eisenstein functions. They do not have poles on $\{\Re s = d/2\}$. We denote $E(s)=(E_1(s), \dots, E_\kappa(s))$. The scattering matrix is defined in the following way. Let $f(s,y)$ be the square matrix whose $i$ line is the zero Fourier coefficient in cusp $Z_i$, at a fixed height $y> a_i$, of $E$. Then
\begin{equation}
f(s,y) = y^s I + \phi(s) y^{1-s}
\end{equation}
Then $\phi(s)$ is the \emph{scattering matrix}, and $\varphi(s)$, its determinant, is the \emph{scattering determinant}. When $\Re s = d/2$, $\phi$ is unitary, and $|\varphi| = 1$.\\

The main result of this note is
\begin{theorem}\label{theorem:QE}
Let $M$ be a cusp manifold whose geodesic flow is \emph{ergodic}. Then, for all compactly supported symbol $\sigma$, we have the following convergence, as $h\to 0$.
\begin{equation}\label{eq:QE}
\begin{split}
\frac{h^d}{4\pi} \int_{\R} &\left| \left\langle \Op(\sigma) E\left(\frac{d}{2} + i\frac{\lambda}{h} \right), E\left(\frac{d}{2} + i\frac{\lambda}{h}  \right) \right\rangle + \frac{\varphi'}{\varphi}\left(\frac{d}{2} + i\frac{\lambda}{h}  \right) \int_{\lambda S^\ast M} \sigma \right| d\lambda \\
&{}+ h^{d+1} \sum_\lambda \left| \langle \Op(\sigma) u_\lambda, u_\lambda \rangle - \int_{h\lambda S^\ast M} \sigma \right| \to 0 .
\end{split}
\end{equation}
\end{theorem}

Observe that contrary to the usual result of Quantum Ergodicity, we do not obtain a \emph{variance}, but a \emph{mean absolute variation}. In some sense, it is because the Eisenstein functions are not $L^2$ that we only get this. The quantization $\Op$ that appears in the theorem, can be any reasonable one, since we use compactly supported symbols. Especially convenient quantizations are the hyperbolic
pseudo-differential ones of \cite{Zelditch-86} that were used in  \cite{Zelditch-91}  and the quantization defined in  \cite{Bonthonneau-2}.

To keep this article brief, we do not review the many constructions and definitions
in \cite{Zelditch-91} and \cite{Bonthonneau-2} but refer the reader to those articles
for further background.

\section{Overview of the proof}

For convenience, we will use the calculus of Pseudo-differential operators of \cite{Bonthonneau-2}. In that article, a quantization procedure $\Op$ was given, that enables one to quantize an appropriate class of symbols. The class of symbols contains in particular the compactly supported functions of $(x,\xi)\in T^\ast M$, the constants, and the functions of the energy $f(|\xi|^2)$, with $f\in C^\infty_c(\R)$.

There are two parts of the proof. First, one has to prove the lemma
\begin{lemma}
The result of the theorem holds if the symbol $\sigma$ has mean value zero in every energy level $\{ |\xi| = $ constant $\}$.
\end{lemma}
As the symbol has mean value zero, and has compact support, the proof of this lemma is very similar to the proof of the compact case \cite{CdV-85,Zelditch-87}. The proof
is given in detail in \cite{Zelditch-91}  for hyperbolic cusp surfaces and the proof 
is essentially the same for all cusp manifolds. Hence we omit the proof.  \\

Let us explain now why the lemma is useful. Let $\sigma\in C^\infty_c(T^\ast M)$. Then, one can define
\begin{equation*}
\overline{\sigma}(\lambda) = \int_{|\xi | = \lambda} \sigma d\mathscr{L}_\lambda
\end{equation*}
Also let $\varepsilon(\sigma, h)$ the quantity in the LHS of \eqref{eq:QE}.

If the average of $\chi \in C^\infty_c(M)$ is $1$ over $M$, then $\sigma$ and $\tilde{\sigma}: (x,\xi) \mapsto \chi(x) \overline{\sigma}(|\xi|)$ have the same average in each energy layer. Hence, to prove the theorem, separate $\sigma$ into $\tilde{\sigma}$ and $\sigma - \tilde{\sigma}$. The latter can be treated using the lemma. Then we find, using the triangular inequality
\begin{equation}
\limsup_{h\to 0} \varepsilon(\sigma, h) = \limsup_{h \to 0} \varepsilon(\tilde{\sigma},h).
\end{equation}

Now, $\tilde{\sigma}$ has a very convenient structure, even more so if we pick carefully the base profile $\chi$, as we will see in the next section. 

However, before we get to the next step, we recall some facts that will be very useful. Let $\Pi^\ast_y$ be the projector on functions whose zero-Fourier mode vanishes for $y_M > y$. We start with the famous Maass-Selberg relation (\cite[section 2]{Sarnak-81}, or \cite{Bonthonneau-4}):
\begin{equation}\label{eq:Maass-Selberg}
\int_{y_M \leq y} \left| E\left( \frac{d}{2} + it \right) \right|^2 = 2 \kappa \log y - \frac{\varphi'}{\varphi} + \Tr \frac{ y^{2i\lambda /h} \phi^\ast - y^{-2i\lambda /h} \phi}{2 i \lambda/h} + \int_{y_M > y} |\Pi^\ast_y E|^2.
\end{equation}
We use the following Poisson formula \cite{Muller-92}:
\begin{equation}\label{eq:Poisson-Formula}
\frac{\varphi'}{\varphi}\left(\frac{d}{2} + it\right) = Q(t) + \sum_{\rho \text{ pole of }\varphi} \frac{2\Re \rho - d}{(\Re \rho - d/2)^2 + (t- \Im \rho)^2}
\end{equation}
for some polynomial $Q$ of power less or equal to $2\lfloor d/2 \rfloor$. Notice all but a finite number of poles of $\varphi$ are on the left of $\{ \Re s < d/2 \}$. We need the Weyl upper bound \cite{Muller-86}:
\begin{equation}\label{eq:Weyl-bound}
\int_{-\lambda}^{\lambda} \frac{\varphi'}{\varphi}\left( \frac{d}{2} + it \right) dt = \mathcal{O}(h^{-d-1}).
\end{equation}
Lastly, if $\sigma$ in the appropriate class of symbols (see \cite{Bonthonneau-2}), 
\begin{equation}\label{eq:Trace-Formula}
\Tr \Pi^\ast_0\Op(\sigma) = h^{-d-1}\int_{T^\ast M} a + \mathcal{O}(h^{-d}|\log h|).
\end{equation}

\section{The case of non-zero mean-value.}

Now, we concentrate on $\tilde{\sigma}(x,\xi) = \chi(x) \overline{\sigma}(|\xi|)$. We can stop considering altogether the $L^2$ eigenfunctions. Indeed, in their contribution to \eqref{eq:QE}, we can replace $\overline{\sigma}(h\lambda_j)$ by $\overline{\sigma}(h \lambda_j) \| u_j \|^2$, and we are back to the case $\overline{\sigma}=0$. Also let $\lambda_0$ be such that $\sigma$ is supported for $|\xi| \leq \lambda_0$.

Now, by the localization of eigenfunctions, 
\begin{equation*}
\begin{split}
\left\langle \Op(\tilde{\sigma})E\left(\frac{d}{2} + i\frac{\lambda}{h}\right),\right. & \left. E\left(\frac{d}{2} + i\frac{\lambda}{h}\right) \right\rangle = \\
	& \int_M \chi \left|E\left(\frac{d}{2} + i\frac{\lambda}{h}\right)\right|^2 \overline{\sigma}\left(\frac{h^2 d^2}{4} +  \lambda^2\right) + \mathcal{O}(h) \|  E \|_{\{\chi \neq 0\}}^2
\end{split}
\end{equation*}
Call $R_1$ the contribution of the remainder in the RHS to $\varepsilon (\tilde{\sigma}, h ) $. Let us first deal with the other term. It will contribute in $\varepsilon (\tilde{\sigma},h)$ by
\begin{equation*}
h^d\int_{-\lambda_0}^{\lambda_0} \left|\overline{\sigma}\left(\frac{h^2 d^2}{4} +  \lambda^2\right)\right| \times \left| \int_M \chi \left|E\left(\frac{d}{2} + i\frac{\lambda}{h}\right)\right|^2 + \frac{\varphi'}{\varphi}\left( \frac{d}{2} + i\frac{\lambda}{h} \right) \right| d\lambda
\end{equation*}
Now, we choose the shape $\chi$ in the following way. We suppose that $\chi$ is constant in $M_0$, and that in the cusps, it writes as $c_\nu \chi_0(y\nu)$, where $\chi_0 \in C^\infty_c(\R)$ equals $1$ close to zero, $\nu$ is a small parameter and $c_\nu \in \R^+$ is a normalization constant. It is chosen so that $\chi$ has mean value $1$, i.e $\int_M \chi = \vol(M)$. Write
\begin{equation*}
\int_M \chi \left|E\left(\frac{d}{2} + i\frac{\lambda}{h}\right)\right|^2 = -c_\nu\int \nu \chi_0'(y\nu) \left\{\int_{y_M\leq y} | E |^2\right\} dy.
\end{equation*}
For $\nu$ small enough, this is well defined, and we can use the Maass-Selberg formula \eqref{eq:Maass-Selberg}:
\begin{equation*}
\begin{split}
\int_M \chi \left|E\left(\frac{d}{2} + i\frac{\lambda}{h}\right)\right|^2 = -c_\nu\int \nu \chi_0'(y \nu) &\left\{2\kappa \log y -  \frac{\varphi'}{\varphi} + \Tr \frac{y^{2i\lambda/h}\phi^\ast - y^{-2i \lambda/h}\phi}{2 i \lambda/h} \right\}dy\\
		&+ \int_M (1-\chi) \left| \Pi^\ast_0 E \right|^2
\end{split}
\end{equation*}
In the RHS, the third term is highly oscillating. Il will only contribute for $\mathcal{O}(h^\infty)$ to the final result, by non-stationary phase. The second will contribute by $-c_\nu \varphi'/\varphi$. The fourth is an integral supported in the cusps. We are left to prove that
\begin{equation*}
\begin{split}
h^d \int_{-\lambda_0}^{\lambda_0} \left|\overline{\sigma}\left(\frac{h^2 d^2}{4} +  \lambda^2\right)\right| \times &\left\{\left| 2 \kappa \int c_\nu \nu \chi_0'(\nu y) . \log y dy\right| + \left|(1-c_\nu)\frac{\varphi'}{\varphi}\left(\frac{d}{2} + i \frac{\lambda}{h} \right) \right| \right\} d\lambda \\
	&+ h^{d+1}\Tr\left\{\Pi^\ast_0 (1-\chi) |\overline{\sigma}|(-h^2\Delta) \right\}.
\end{split}
\end{equation*}
goes to $0$ with $h$. A direct computations shows this for the first term. Let us call the third term $R_2$. The second term gives
\begin{equation*}
h^d|1 - c_\nu| \int_{-\lambda_0}^{\lambda_0} h^d \overline{\sigma}\left(\frac{h^2 d^2}{4} +  \lambda^2\right)\left|\frac{\varphi'}{\varphi}\left(\frac{d}{2} + i \frac{\lambda}{h} \right)\right| 
\end{equation*}
In the decomposition \eqref{eq:Poisson-Formula} the sum over the resonances is a function that is always negative for big values of $t$, and the polynomial part is explicit, so we know that
\begin{equation*}
\int_{-\lambda_0}^{\lambda_0} \left|\frac{\varphi'}{\varphi}\left(\frac{d}{2} + i \frac{\lambda}{h} \right)\right| + \frac{\varphi'}{\varphi}\left(\frac{d}{2} + i \frac{\lambda}{h} \right) d\lambda = \mathcal{O}(h^{-d}).
\end{equation*}
Combining this with the Weyl bound \eqref{eq:Weyl-bound}, we arrive at
\begin{equation*}
|1 - c_\nu| \int_{-\lambda_0}^{\lambda_0} h^d \overline{\sigma}\left(\frac{h^2 d^2}{4} +  \lambda^2\right)\left|\frac{\varphi'}{\varphi}\left(\frac{d}{2} + i \frac{\lambda}{h} \right)\right|  = |1 - c_\nu|\mathcal{O}(1).
\end{equation*}

Now, we come back to $R_1$ and $R_2$. First, $R_1$ can be bounded by some $\mathcal{O}(h)\Tr \Op(\sigma')$, with $\sigma'$ compactly supported (on a set slightly bigger than $\tilde{\sigma}$. Then, we can apply formula \eqref{eq:Trace-Formula}. For $R_2$, we can \emph{also} apply the formula \eqref{eq:Trace-Formula} (the symbol also is in the quantizable class). Hence, we find that $R_1$ and $R_2$ contribute to the $\limsup$ by 
\begin{equation*}
\mathcal{O}(1)\vol\{ y_M \geq 1/\nu \}.
\end{equation*}
At last, let us observe that when $\nu \to 0$, the assumption that $\chi$ has mean value $1$ implies that $c_\nu \to 1$. Actually, this can be made more precise, as $1-c_\nu = \mathcal{O}(\vol\{ y > 1/\nu \})$. We deduce that for some constant $C >0$,
\begin{equation*}
\limsup_{h\to 0}\varepsilon(\tilde{\sigma}, h)  \leq C \vol\{ y \geq 1/\nu \}.
\end{equation*}
As this does not depend on $\nu$, we can take $\nu \to 0$, and this ends the proof.

\end{document}